\documentclass{amsart}

\usepackage{amssymb,amsmath}
\usepackage[cmtip,all]{xy}
\usepackage{enumitem}

\newtheorem{theorem}{Theorem}[section]
\newtheorem{lemma}[theorem]{Lemma}

\theoremstyle{definition}
\newtheorem{definition}[theorem]{Definition}

\theoremstyle{remark}

\numberwithin{equation}{section}

\newcommand{\dcont}{d_{{\rm c}}}

\newcommand{\R}{\mathbb{R}}

\newcommand{\Z}{\mathbb{Z}}

\newcommand{\osc}{{\rm osc}}
\newcommand{\id}{\textup{id}}

\begin{document}

\title[Topological strictly contact isotopies]{Uniqueness of contact Hamiltonians of topological strictly contact isotopies}

\author[A.~Banyaga]{Augustin Banyaga}

\address{Penn State University, University Park, USA}
\email{banyaga@math.psu.edu}

%    author two information
\author[P.~Spaeth]{Peter Spaeth}

\address{GE Global Research, Niskayuna, NY, USA}
\email{peter.spaeth@ge.com}

%\subjclass[2000]{Primary }
%    The 2010 edition of the Mathematics Subject Classification is
%    now available.  If you are citing a classification from the
%    new scheme, use the following input coding instead.

\date{\today}

\begin{abstract}
We prove that for regular contact forms there exists a bijective correspondence between the $C^0$ limits of sequences of smooth strictly contact isotopies and the limits with respect to the contact distance of their corresponding Hamiltonians.
%The uniqueness of the topological isotopy of a topological contact Hamiltonian implies the $C^0$ rigidity of strictly contact diffeomorphisms.
\end{abstract}

\maketitle


\section{Introduction}
Viterbo \cite{viterbo:ot06} proved that if two sequences $\Phi_{H^{j}}$ and $\Phi_{F^j}$ of Hamiltonian isotopies converge in the $C^0$-topology to the same limit (a continuous Hamiltonian system), and the corresponding sequences $H^j$ and $F^j$ of normalized Hamiltonians converge uniformly to continuous functions $H$ and $F$, then $H = F$.
Buhovsky-Seyfaddini subsequently \cite{buhovsky:ugh10} found a generalization of Viterbo's result where it is sufficient for the sequences $H^j$ and $F^j$ of normalized Hamiltonians to converge in the Hofer norm.

The main goal of this note is to prove the uniqueness of the topological contact Hamiltonian of a topological strictly contact isotopy of a compact regular contact manifold.

\begin{theorem}\label{thm:uniqueness}
If two sequences $\Phi_{H^{j}}$ and $\Phi_{F^j}$ of smooth strictly contact isotopies of a compact manifold $M$ with a regular contact form $\alpha$ uniformly converge as $j \to \infty$ to a continuous isotopy $\Phi = \{ \phi_t\}$ of homeomorphisms of $M$, and the corresponding sequences $H^{j}$ and $F^j$ of generating Hamiltonians converge to functions $H$ and $F:M \times [0,1] \to \R$ with respect to the contact norm, then $H = F$.
\end{theorem}

Here the contact norm \cite{banyaga-donato} of a function $H:M\times [0,1] \to \R$ is given by
\begin{align}\label{eqn:contact-norm}
 \| H \| = \int_0^1 \left( \max_{x \in M} H(x,t) - \min_{x \in M} H(x,t) + \frac{ \left| \int_M H(x,t) \alpha \wedge (d\alpha)^m \right| }{\int_M \alpha \wedge (d\alpha)^m} \right) dt.
 \end{align}

Theorem~\ref{thm:uniqueness} also has useful applications in hydrodynamics (helicity and asymptotic Hopf invariant \cite{mueller:hvf11}).

\section{Review of contact geometry}
Let $M$ be a smooth manifold and ${\rm Diff}(M)$ denote the group of diffeomorphisms of $M$ with the compact-open topology.
We write ${\rm Diff}_0(M)$ for its identity component.
A {\bf smooth isotopy} of $M$ is a smooth map $\Phi: [0,1] \times M \to M$ such that for each time $t \in [0,1]$, the map 
	$ \Phi(t, \cdot) = \phi_t :M \to M$
is a diffeomorphism and $\phi_0 = \id$.
We denote the group of smooth isotopies of $M$ by $P{\rm Diff}(M,\id)$.
Every smooth isotopy $\Phi = \{ \phi_t \}$ is generated by the family of vector fields
\begin{eqnarray}\label{eq:vector-field}
\dot \phi_t = \frac{d\phi_t}{dt} \circ (\phi_t)^{-1}.
\end{eqnarray}

A {\bf contact form} $\alpha$ on a smooth manifold $M$ of dimension $2m+1$ is a 1-form such that
	$ \alpha \wedge (d\alpha)^m \neq 0$.
Associated to $\alpha$ is a unique vector field $R$, called the {\bf Reeb vector field}, such that
	$ \iota(R)\alpha = 1 \ {\rm and} \ \iota(R)d\alpha = 0$,
where $\iota(X)\beta$ denotes the contraction of a vector field $X$ with a differential form $\beta$.
A smooth isotopy $\{ \phi_t \}$ is a {\bf smooth strictly contact isotopy} if
\begin{eqnarray}
\phi_t^*\alpha = \alpha
\end{eqnarray}
for all $t \in [0,1]$.
This means that $\dot \phi_t$ satisfies $L_{\dot \phi_t}\alpha = 0$, where $L_X$ denotes the Lie derivative along $X$.
A vector field $X$ such that $L_X \alpha = 0$ is called a {\bf strictly contact vector field}.
The Lie algebra $\mathcal {L}(M,\alpha)$ of strictly contact vector fields is isomorphic to the set 
	\[ C^\infty_b(M) = \{ f:M \to \R \mid R \cdot f = df(R) =  0 \}\]
of smooth {\bf basic functions} on $M$.

The isomorphism $\mathcal{L}(M,\alpha) \to C^\infty_b(M)$ sends $X$ to the function $\iota(X)\alpha$, while the inverse sends $H$ to the vector field $X_H$ defined by the two equations $\iota(X_H)d\alpha = -dH$ and $\iota(X_H)\alpha = H$.

If $\Phi$ is a smooth strictly contact isotopy, then the family of functions $H = \{ H_t \} =\{ \iota(\dot \phi_t)\alpha \}$ is called the {\bf contact Hamiltonian} of $\Phi$, and $\dot \phi_t = X_{H_t}$.
We say that $H$ {\bf generates} $\Phi$, and write $\Phi = \Phi_H = \{\phi^t_H\}$.

The composition $\Phi_H \circ \Phi_F = \{ \phi^t_H \circ \phi^t_F \}$ is generated by the Hamiltonian
\begin{eqnarray}\label{eq:hamil-prod}
H\#F (t,x) = H(t,x) + F\left(t, (\phi^t_H)^{-1}(x)\right),
\end{eqnarray}
and the inverse isotopy $\Phi_F^{-1} = \{ (\phi_H^t)^{-1} \}$ is generated by the Hamiltonian 
\begin{eqnarray}\label{eq:hamil-inv}
\overline F(t,x) = -F(t, \phi^t_F(x)).
\end{eqnarray}

The group of smooth strictly contact isotopies of $(M,\alpha)$ is denoted $P{\rm Diff}(M,\alpha)$ and
	\[ {\rm Diff}(M,\alpha) = \{ \phi \in {\rm Diff}(M) \mid \phi^*\alpha = \alpha\}\]
denotes the group of strictly contact diffeomorphisms of $(M,\alpha)$.
Its identity component ${\rm Diff}_0(M,\alpha)$ consists of the time-one maps of strictly contact isotopies.

We denote the set of pairs $(\Phi_H, H)$, where $\Phi_H$ is a smooth strictly contact isotopy generated by the smooth basic contact Hamiltonian $H$ by $\mathcal{C}(M,\alpha)$. 
It is a topological group under the operations
\begin{align*}
	&(\Phi, H)\circ (\Psi, F) = (\Phi \circ \Psi, H \#F) = (\{ \phi^t_H \circ \psi^t_F\},  H_t+ F_t \circ (\phi_t)^{-1}),\ \textup{and} &\\
	&(\Phi, H)^{-1} = (\Phi^{-1}, \overline H) =  (\{ (\phi^t_H)^{-1} \}, -H_t\circ \phi_t).&
\end{align*}

If $\varphi \in {\rm Diff}(M,\alpha)$ then the isotopy $\varphi^*\Phi_H = \{ \varphi^{-1}\circ \phi^t_H \circ \varphi \}$ is generated by the Hamiltonian $H^\varphi$ defined by
\begin{eqnarray}\label{eqn:auto-conjugate}
H^\varphi_t = H_t \circ\varphi .
\end{eqnarray} 

The isomorphism $\mathcal{L}(M,\alpha) \cong C^\infty_b(M)$ allows us to define a norm on $\mathcal{L}(M,\alpha)$ and a distance on $P{\rm Diff}(M,\alpha)$ as follows.
The 1-form $\alpha$ defines a canonical volume form $\nu_{\alpha} = \alpha \wedge (d\alpha)^{m}$ on $M$. 
For a basic function $f\in C^\infty_b(M)$, set
\begin{eqnarray}\label{eq:norm}
\|f\|_{\rm c} = \osc(f) + |c(f)|,
\end{eqnarray}
where $\osc (f) = \max_{x \in M} f(x) - \min_{x \in M} f(x)$,
and $c(f) = \frac{1}{\int_M  \nu_\alpha}\int_M f\,  \nu_\alpha$ denotes the average value of the function $f$.

Equation (\ref{eq:norm}) then gives rise to a distance on $P{\rm Diff}(M,\alpha)$.
For strictly contact isotopies $\Phi_H$ and $\Phi_F$, define
\begin{eqnarray}\label{eq:length}
D(\Phi_H, \Phi_F) := \|F - H\| = \int_0^1 \| F_t - H_t \|_{\rm c}\, dt.
\end{eqnarray}

Finally a contact form $\alpha$ is called {\bf regular} if its Reeb field generates a free $S^1$-action on $M$.
The quotient manifold $B = M/S^1$ is an integral symplectic manifold $(B, \omega)$, and the canonical projection $\pi : M \to B$ is a principal $S^1$-fibration with $\pi^*\omega = d\alpha$.
The contact form $\alpha$ is a connection in the principal bundle $\pi$, and $\pi^*\omega$ is its curvature.
Moreover, the cohomology class $[\omega]$ of $\omega$ is integral.
Conversely, any closed symplectic manifold $(B,\omega)$ such that $[\omega]$ is an integral cohomology class is the base of a principal circle bundle $\pi : M \to B$ whose Chern class is $[\omega]$; $M$ carries a regular contact form that is a connection form with curvature $\pi^*\omega$ \cite{boothby:oc58}. 
Regular contact forms are of great interest in Mechanics \cite{souriau:sd69}.

\section{Topological strictly contact isotopies}
If $\Phi = \{ \phi_t \}$ and $\Psi= \{ \psi_t \}$ are continuous isotopies of $M$, then we define
	\[ \overline {d}(\Phi,\Psi) = \max_{t\in [0,1]} \max\{d_0(\phi_t,\psi_t), d_0((\phi_t)^{-1},(\psi_t)^{-1})\},\]
where $d_0(\phi,\psi) = \max_{x \in M}d_{g_M}(\phi(x), \psi(x))$ is the distance induced by any Riemannian metric $g_M$ on $M$, and $\phi$ and $\psi \in {\rm Homeo}(M)$ are homeomorphisms.
Then $\overline d$ is a complete metric that induces the $C^0$-topology on the group $P{\rm Homeo}(M,\id)$ of continuous isotopies of homeomorphisms of $M$ based at the identity.
Note that this topology is independent of the choice of the metric $g_M$.

\begin{definition}
The {\bf contact topology} on $P{\rm Diff}(M,\alpha)$ is the metric topology induced by the contact distance $\dcont$, where the distance between two strictly contact isotopies $\Phi_F$ and $\Phi_H$ is
\begin{eqnarray*}
\dcont(\Phi_H, \Phi_F) & := & \overline {d}(\Phi_H,\Phi_F) + D(\Phi_H, \Phi_F)\\
& = & \overline {d}(\Phi_H,\Phi_F) + \|F - H\|\\
& = & \overline {d}(\Phi_H,\Phi_F) + \|\overline H \# F\|.
\end{eqnarray*}
\end{definition}
The third equality follows from (\ref{eq:hamil-prod}), (\ref{eq:hamil-inv}), and the fact that $\phi^t_H$ preserves $\nu_\alpha$.

\begin{definition}
An isotopy $\Phi$ of homeomorphisms of $M$ is called a {\bf topological strictly contact isotopy} if there exists a $\dcont$-Cauchy sequence $\Phi_{H_j}$ of smooth strictly contact isotopies that uniformly converges to $\Phi$.
The set of topological strictly contact isotopies is denoted by $P{\rm Homeo}(M,\alpha)$.

A homeomorphism $\phi$ of $M$ is called a {\bf strictly contact homeomorphism} if it is the time-one map of a topological strictly contact isotopy.
We denote the set of all such homeomorphisms by ${\rm Homeo}(M,\alpha)$.

Let $\Phi_{H_j}$ be a $\dcont$-Cauchy sequence of smooth strictly contact isotopies.
The limit with respect to the norm in equation (\ref{eqn:contact-norm}) of the sequence $H_j$ is called a {\bf topological contact Hamiltonian}.
%Then $\lim_{\| \cdot \|} H_j$ is called a {\bf topological contact Hamiltonian}.
\end{definition}

A topological contact Hamiltonian function $H$ can be regarded as an element of the space $L^1([0,1], C^0(M))$ of $L^1$-functions on the unit interval with values in the space $C^0(M)$ of continuous functions on $M$.
Let $\mathcal{TC}(M,\alpha)$ denote the set of limits $(\Phi, H)$ of $\dcont$-convergent sequences $(\Phi_{H_j}, H_j) \in \mathcal{C}(M,\alpha)$.
The next result shows that non-smooth time reparametrizations of the Reeb flow are elements of $\mathcal{TC}(M,\alpha)$.

\begin{lemma}\label{lem:top_reeb}
Let $(\{ \varphi_R^t \}, {\bf 1}) \in \mathcal{C}(M,\alpha)$ denote the Reeb flow on any regular $(M,\alpha)$, where ${\bf 1}$ denotes the constant function equal to $1$ on $M$.
If $f \in L^1([0,1])$ is any $L^1$ function of time, then the pair $( \varPhi_R^a , f)$ is an element of $\mathcal{TC}(M,\alpha)$, where $\{ \varphi_R^{a(t)} \}$ denotes the reparametrization of the Reeb flow $\{ \varphi_R^t \}$ at time $a(t) = \int_0^t f(\tau) d\tau$. 
\end{lemma}

\begin{proof}
Using the regularity of $\alpha$ we may assume that the $C^0$-metric $\overline d$ on $P{\rm Homeo}(M,\id)$ is induced by a Riemannian metric of the form $g_M = \pi^*g_B + \alpha \otimes \alpha$, where $g_B$ is a Riemannian metric on $B$.

Let $f_k : [0,1] \to \R$ be a sequence of smooth functions that converges to $f$ in $L^1([0,1])$, and denote by $a_k(t) = \int_0^t f_k(\tau) d\tau$.
We claim: (a) that $f_k$ converges to $f$ under the contact norm, and (b) that $\varPhi^{a_k}_R$ converges uniformly to $\varPhi_R^{a}$.

Note that as functions of time $t$ alone, the oscillations over $M$ of $f_k$ and $f$ vanish, hence
\begin{align*}
\|f_k - f \| = \int_0^1 \left( \frac{ \left| \int_M (f_k(t) - f(t)) \alpha \wedge (d\alpha)^m \right| }{\int_M \alpha \wedge (d\alpha)^m} \right) dt
\le \int_0^1 \left| f_k(t) - f(t) \right| dt.
\end{align*}
This proves (a).

To prove (b), we begin by fixing $t \in [0,1]$ and $x \in M$.
Consider the following curve $\sigma: [0,1] \to M$ beginning and ending at the points $\varphi^{a(t)}_R(x)$ and $\varphi^{a_k(t)}_R(x)$, which lie in the same fiber over the point $\pi(x) \in B$. 
	\[ \xymatrix@C=1pc@R=1pc{ \varphi^{a_k(t)}_R(x)\!\!\!\!\!\!\!\!  & \bullet \\  \varphi_R^{a(t)}(x)\!\!\!\!\!\!\!\! & \bullet \ar@/_1pc/[u]_{\sigma(s) = \varphi_R^{(1-s)a(t) + s a_k(t)(x)} }} \]

Because $\sigma$ follows the direction of the fiber, its tangent vector $\sigma'(s) = (a_k(t) - a(t))R_{\left. \right| \sigma(s)}$ satisfies
	\[ |\sigma'(s)|_{g_M}  = \sqrt {g_M(\sigma'(s), \sigma'(s) )} = \sqrt {(a_k(t) - a(t))^2 } = |a_k(t) - a(t)|. \]
Therefore the length of the curve $\sigma$ satisfies
	\[ {\rm length}(\sigma) = \int_0^1 \left| \sigma'(s)\right|_{g_M} ds = \int_0^1 \left| a_k(t) - a(t) \right| ds = \left| a_k(t) - a(t) \right|.\]
Using the definitions of $a_k$ and $a$ gives
\begin{align*}
{\rm length}(\sigma) = & \left| a_k(t) - a(t) \right| = \left| \int_0^t f_k(\tau) - f(\tau) d\tau \right|\\
& \le \int_0^1 \left| f_k(\tau) - f(\tau) \right| d\tau =  \|f_k - f \|_{L^1([0,1])}.
\end{align*}
Since the distance between $\varphi^{a(t)}_R(x)$ and $\varphi^{a_k(t)}_R(x)$ equals the infimum of the lengths of curves joining these points we may conclude that
	\[d_{g_M}\left( \varphi^{a(t)}_R(x), \varphi^{a_k(t)}_R(x) \right) \leq \|f_k - f\|_{L^1([0,1])}.\] 
Taking the maximum over $x \in M$ and $t \in [0,1]$ we conclude that the isotopies satisfy $d(\varPhi^{a_k}_R, \varPhi^{a}_R) \leq \|f_k - f\|_{L^1([0,1])}$.
This proves (b).

To conclude the proof, note that the inverse of the smooth pair $(\{ \varphi_R^{a_k(t)} \} , f_k)$ equals $( \{ \varphi_R^{-a_k(t)} \}, -f_k)$; therefore the proofs that $- f_k$ and $(\varPhi_R^{a_k})^{-1}$ converge to $-f$ and $(\varPhi_R^{a})^{-1}$ respectively are similar.
\end{proof}

Because a strictly contact diffeomorphism $\phi$ preserves the contact form $\alpha$, the techniques of \cite{oh:tg07} can be adapted to prove the next two results.

\begin{theorem}\label{thm:topgroup}
Let $M$ be a compact manifold with a contact form $\alpha$.
The set $\mathcal{TC}(M,\alpha)$ admits the structure of a topological group under the operations
\begin{align*}
&(\Phi, H)\circ (\Psi, F) = (\Phi \circ \Psi, H \#F) = (\{ \phi^t_H \circ \psi^t_F\},  H_t+ F_t \circ (\phi_t)^{-1}),\ and &\\
&(\Phi, H)^{-1} = (\Phi^{-1}, \overline H) =  (\{ (\phi^t_H)^{-1} \}, -H_t\circ \phi_t).&
\end{align*}
$\mathcal{TC}(M,\alpha)$ contains $\mathcal{C}(M,\alpha)$ as a subgroup.
\end{theorem}

Recall that $\Phi$ are $\Psi$ are isotopies of homeomorphisms of $M$ that are $C^0$-limits of sequences of strictly contact isotopies, say $\Phi_j = \Phi_{H_j}$ and $\Psi_j = \Psi_{F_j}$, whose generating Hamiltonians $H_j$ and $F_j$ converge to topological contact Hamiltonians $H$ and $F$.
Since the sequence $\Phi_j^{-1} \circ \Psi_j$ obviously $C^0$-converges to $\Phi^{-1} \circ \Psi$ to prove Theorem \ref{thm:topgroup} we must show that the corresponding sequence of contact Hamiltonians $\overline{H}_j\# F_j$ is Cauchy with respect to the contact norm.
The $t$ variable is suppressed in the next Lemma and the proof of Theorem \ref{thm:topgroup}.
\begin{lemma}\label{lemma:uniform}
Let $\epsilon>0$ be given.
For any integer $j_0$, there exists $N_0$, such that for any $j,k \geq N_0$ we have $\|F_{j_0} \circ \phi_j - F_{j_0}\circ\phi_k\| < \epsilon$.
\end{lemma} 
\begin{proof}
The uniform continuity of $F_{j_0}$ implies that there exists $\delta = \delta(F_{j_0})> 0$, such that $\overline d(\phi_j, \phi_k) < \delta$ implies $\|F_{j_0} \circ \phi_j - F_{j_0}\circ \phi_k\| < \epsilon$.
Such $N_0$ exists since $\overline d(\phi_j, \phi_k) \to 0.$
\end{proof}

\begin{proof}[Proof of Theorem~\ref{thm:topgroup}]
Let $\epsilon> 0$ be given.
There exist $j_1$ and $j_2$ such that 
	\[ j,k \ge j_1 \implies \|H_j - H_k\|  \le  \epsilon\ \textup{and}, j,k \ge j_2 \implies \|F_j - F_k\|  \le  \epsilon. \]
Let $j_0 \geq \max\{j_1,j_2\}$, and $N_0$ be large enough to satisfy the preceding lemma applied to either function $F_{j_0}$ or $H_{j_0}$ with the given $\epsilon$.
Let $j,k \geq N_0$.
Then,
	\[ \|\overline{H}_j\#F_j - \overline{H}_k\# F_k\| = \|(F_j - H_j)\circ \phi_j - (F_k-H_k)\circ \phi_k\| \]
by equations (\ref{eq:hamil-prod}) and (\ref{eq:hamil-inv}).
Adding and subtracting the quantities $F_k\circ \phi_j$ and $H_j \circ \phi_k$, and applying the triangle inequality yields
\begin{eqnarray*}
\|\overline{H}_j\# F_j - \overline{H}_k\#F_k\| & \leq &  \|F_j\circ \phi_j - F_k\circ \phi_j\|  + \|F_k\circ \phi_j - F_k\circ \phi_k\| \\
& & + \|H_k\circ \phi_k - H_j \circ \phi_k\| + \|H_j\circ \phi_k - H_j\circ \phi_j\|.
\end{eqnarray*}
Because $\phi_j$ and $\phi_k$ both preserve $\alpha$, we have 
\[ \|F_j\circ \phi_j - F_k\circ \phi_j\| = \|F_j - F_k\|< \epsilon,\]
and
\[ \|H_k\circ \phi_k - H_j\circ \phi_k\| = \|H_k - H_j\|<\epsilon.\]

Now consider the term $\|F_k \circ \phi_j - F_k \circ \phi_k\|$.
By the triangle inequality,
\begin{eqnarray*}
\|F_k\circ \phi_j - F_k \circ \phi_k\| & \leq &  \|F_k \circ \phi_j - F_{j_0}\circ \phi_j\| + \|F_{j_0}\circ \phi_j - F_{j_0}\circ \phi_k\|\\
& & +  \|F_{j_0}\circ \phi_k - F_k\circ \phi_k\|\\
& = &  \|F_k - F_{j_0}\| + \|F_{j_0}\circ \phi_j - F_{j_0}\circ \phi_k\| + \|F_{j_0} - F_k\| \\
& <  &  3\, \epsilon.
\end{eqnarray*} 
Similarly one obtains the inequality $\|H_j \circ \phi_k - H_j \circ \phi_j\| < 3\, \epsilon.$
Putting this all together, we have shown that $\|\overline{H}_j\#F_j - \overline{H}_k\#F_k\| < 8\, \epsilon$.

We now show that composition is continuous with respect to the topology induced from the obvious extension of the contact metric $\dcont$ on the set $\mathcal{TC}(M,\alpha)$.
Let $\left( \Phi_j, H_j\right), \left( \Psi_j,F_j\right)\in \mathcal{TC}(M,\alpha)$ be sequences such that,
\begin{eqnarray*}
&(\Phi_j,H_j) \longrightarrow (\Phi,H) \ {\rm and}\ (\Psi_j,F_j) \longrightarrow (\Psi,F)
\end{eqnarray*}
as $j \to \infty$, for some $(\Phi, H)$ and $(\Psi, F) \in \mathcal{TC}(M,\alpha)$, where convergence is with respect to the metric $\dcont$.

Again, $\Phi_j \circ \Psi_j$ converges to $\Phi \circ \Psi$ with respect to $\overline d$.
Moreover
\begin{eqnarray*}
\|H_j + F_j \circ \phi^{-1} - (H + F \circ \phi^{-1})\| \leq \|H - H_j\| + \| F_j \circ \phi_j^{-1} - F \circ \phi^{-1},\| 
\end{eqnarray*}
and this in turn can be estimated using the triangle inequality as
\begin{eqnarray*}
\|H - H_j\| + \| F_j \circ \phi_j^{-1} - F \circ \phi^{-1}\| \leq \|H - H_j\|
+ \|F_j - F\|  + \|F \circ \phi_j^{-1} - F  \circ \phi^{-1}\|.
\end{eqnarray*}
Thus
\begin{eqnarray}\label{eqn:product_estimate}
\|H_j + F_j \circ \phi^{-1} - (H + F \circ \phi^{-1})\| & \leq & \|H - H_j\| + \|F - F_j\|\\ \nonumber
& + & \|F \circ \phi_j^{-1} - F  \circ \phi^{-1}\|.
\end{eqnarray}

Both $\|H - H_j\|$ and $\|F - F_j\|$ approach zero by assumption.
By the triangle inequality the third term above satisfies
 \begin{eqnarray*}
 \|F \circ \phi_j^{-1} - F  \circ \phi^{-1}\|  & \leq & \| F \circ \phi_j^{-1} - F_j \circ \phi_j^{-1}\| \\
 & & +  \| F_j \circ \phi_j^{-1} - F_j \circ \phi^{-1}\| + \|(F_j - F) \circ \phi^{-1} \|.
 \end{eqnarray*}
The first and third summands approach zero by assumption.
By Lemma~\ref{lemma:uniform}, $\| F_j \circ \phi_j^{-1} - F_j \circ \phi^{-1}\|  \to 0$ as $j \to \infty$.
This proves that under the $\dcont$-metric we have
	\[ (\phi_j\circ \psi_j, H_j\# F_j) \to  (\phi \circ \psi, H + F \circ \phi^{-1}), \]
and so the composition $\mathcal{TC}(M,\alpha) \times \mathcal{TC}(M,\alpha)\to \mathcal{TC}(M,\alpha)$,	
	\[ ((\Phi,H), (\Psi, F)) \mapsto (\Phi \circ \Psi, (H \# F) = H + F \circ (\phi)^{-1}) \]
is continuous.
The proof that the inverse is continuous is similar. 
\end{proof}

As an extension of the case of smooth isotopies and contact Hamiltonians we have the following.

\begin{theorem}\label{thm:group}
The group ${\rm Homeo}(M,\alpha)$ is a normal subgroup of $\overline {\rm Diff}(M,\alpha)$, where
	\[ \overline {\rm Diff}(M,\alpha) \subset {\rm Homeo}(M)\]
denotes the closure of ${\rm Diff}(M,\alpha)$ under the $C^0$-topology.

More precisely, if $\psi \in \overline {\rm Diff}(M,\alpha)$ and $(\Phi, H) \in \mathcal{TC}(M,\alpha)$, then $(\psi^*\Phi, H^\psi)\in \mathcal{TC}(M,\alpha)$, where $ \psi^*\Phi = \{ \psi^{-1} \circ \phi_t \circ \psi \}$ and $\ H^\psi_t = H_t \circ \psi$.
\end{theorem}

\begin{proof}
The definitions imply that $\Phi = \{\phi_t\}$ is the $C^0$-limit of a sequence of smooth strictly contact isotopies $\Phi_{H_i}$ such that $H_i$ is a $\|\cdot \|$-Cauchy sequence of smooth basic functions, and $\psi$ is the $C^0$-limit of a sequence $\psi_i$ of strictly contact diffeomorphisms of $(M,\alpha)$.

Because $\psi_i$ and $\Phi_{H_i}$ both $C^0$-converge to $\psi$ and $\Phi$ respectively, the composition $\psi_i^{-1}\circ \Phi_{H_i} \circ \psi_i$ converges in $C^0$ to $\psi^{-1} \circ \Phi \circ \psi$.
%, and in particular \[ \psi^{-1} \circ \phi \circ \psi = \psi^{-1} \circ \phi_1 \circ \psi. \]

The smooth basic contact Hamiltonian $\psi_i^* H_i$ generates the smooth strictly contact isotopy $\psi_i^*\Phi_{H_i} = \psi_i^{-1} \circ \Phi_{H_i}\circ \psi_i$ and we claim that the sequence of basic Hamiltonians $\psi_i^*H_i$ is Cauchy with respect to the contact norm $\| \cdot \|$.

Consider $\| \psi_i^*H_i - \psi_j^*H_j \|$.
By the triangle inequality and the fact that each $\psi_i$ preserves the contact form $\alpha$,
\begin{eqnarray*}
\| \psi_i^*H_i - \psi_j^*H_j \| & \leq & \| \psi_i^*H_i - \psi_j^*H_i \| + \| \psi_j^*H_i - \psi_j^*H_j\|\\
& = & \| H_i - (\psi_j^{-1}\circ \psi_i)^*H_i\| + \| H_i - H_j\|.
\end{eqnarray*}

The argument that concludes the proof of Theorem~\ref{thm:topgroup} implies that the first term approaches zero by Lemma~\ref{lemma:uniform} and the triangle inequality, while the second term approaches zero by assumption.
Hence $\psi^*\Phi = \psi^{-1} \circ \Phi \circ \psi$ is a topological strictly contact isotopy, and ${\rm Homeo}(M,\alpha)$ is a normal subgroup of $\overline {\rm Diff}(M,\alpha)$.
\end{proof}

\begin{theorem}
The group ${\rm Homeo}(M,\alpha)$ is contained in the identity component of the group of $\overline{\nu}_\alpha$-preserving homeomorphisms of $M$, where $\bar \nu_\alpha$ is the measure on $M$ induced by the volume form $\nu_\alpha = \alpha \wedge (d\alpha)^m$.
\end{theorem}

\begin{proof}
Let $\phi \in {\rm Homeo}(M,\alpha)$, $U \subset M$ be a measurable subset, and $\chi_U$ be the characteristic function of $U$.
Let $\phi_k$ be a sequence of strictly contact diffeomorphisms that $C^0$-converges to $\phi$.
For all $k$,
\begin{align*}
|U| = \int_M \chi_U\cdot  \nu_\alpha = \int_M \phi_k^* (\chi_U\cdot  \nu_\alpha)
\end{align*}
by the change of variables formula, and
\begin{align*}
\int_M \phi_k^* (\chi_U\cdot  \nu_\alpha)=\int_M (\chi_U \circ \phi_k)\cdot \nu_\alpha
\end{align*}
since $\phi_k^* \alpha = \alpha$.

Therefore
	\[ |U| = \lim_{k \to \infty} \int_M (\chi_U \circ \phi_k) \cdot  \nu_\alpha.\]
By Fatou's lemma
\begin{align*}
|U| =\lim_k \int_M (\chi_U \circ \phi_k) \nu_\alpha
\leq \int_M \lim_{k\to \infty} \chi_U \circ \phi_k \ \nu_\alpha
= \int_M (\chi_U \circ \phi) \nu_\alpha= | \phi(U)|.
\end{align*}
Hence $|U| \leq |\phi(U)|$ for all measurable subsets $U \subset M$.
%Similarly $|U | \leq |\phi^{-1}(U)|$.  

Repeating the same argument with $U$ replaced by $V = \phi(U)$ yields
	\[ |\phi(U)| = | V | \leq |\phi^{-1}(V)| = |U|. \]
Hence $|\phi(U)| = |U|$.
\end{proof}

\section{The proof of Theorem~\ref{thm:uniqueness}}\label{sec:uniqueness}

The following lemma is an important ingredient used to prove Theorem \ref{thm:uniqueness}.
Its proof combines arguments of Hofer-Zehnder \cite{hofer:si05} and Oh-M\"uller \cite{oh:tg07} with the contact energy-capacity inequality from \cite{mueller:tcd11}.
Note that $\alpha$ need not be regular.

\begin{lemma}\label{thm:non-triviality}
If two sequences $\Phi_{H^{j}}$ and $\Phi_{F^{j}}$ of strictly contact isotopies of a compact manifold $M$ with a contact form $\alpha$ uniformly converge as $j \to \infty$ to two isotopies $\Phi = \{\phi_t\}$ and $\Psi = \{ \psi_t \}$ of homeomorphisms of $M$, and  the corresponding sequences of basic Hamiltonians $H^{j}$ and $F^j$ converge to the same function $H: M \times [0,1] \to \R$, then $\Phi = \Psi$, i.e. $\phi_t = \psi_t$ for all $t \in [0,1]$.
\end{lemma}

\begin{proof}%[Proof of Lemma~\ref{thm:non-triviality}]
It suffices to show that if $H=0$, then $\Phi$ is the constant isotopy at the identity.
We argue by contradiction and assume that $\Phi \neq \textrm{Id}$.
By a standard reparametrization argument, without loss of generality we may assume that the time-one map $\varphi = \phi_1$ of the continuous isotopy $\Phi$ is not equal to the identity.
Hence there exists a point $x \in M$ such that $\varphi(x) \neq x$.
By continuity, there exists a neighborhood $U$ containing $x$ such that $\varphi(U) \cap U = \emptyset$.
The $C^0$ convergence of the sequence $\Phi_{H_j}$ to $\Phi$ implies that for all sufficiently large $j$ we have
	\[ \phi^1_{H_j} (U) \cap U = \emptyset. \]
The contact energy-capacity inequality from \cite[Theorem~1.1]{mueller:tcd11} implies that there exists $C > 0$ independent of $j$ such that
	\[ 0 < C \le \| H_j \|,\]
which leads to a contradiction, because $\| H_j \| \to 0$ as $j \to \infty$ by assumption.
\end{proof}

\begin{proof}[Proof of Theorem~\ref{thm:uniqueness}]
Recall that we have sequences $\Phi_{H^j}$ and $\Phi_{F^j}$ of smooth strictly contact isotopies of $(M,\alpha)$ that $C^0$ converge as $j \to \infty$ to the same continuous isotopy of homeomorphisms $\Phi$, and such that the corresponding generating contact Hamiltonians $H^j$ and $F^j$ converge with respect to the contact norm (\ref{eqn:contact-norm}) to topological Hamiltonians $H$ and $F \in L^1([0,1], C^0(M))$.
We must establish that $H = F$.

Consider the prequantization bundle $S^1 \to M \stackrel{\pi}{\to} B$.
We again assume that the $C^0$-metric $\overline d$ on $P{\rm Homeo}(M,\id)$ is induced by any Riemannian metric of the form $g_M = \pi^*g_B + \alpha \otimes \alpha$, where $g_B$ is a Riemannian metric on $B$.

Since $\Phi$ is the limit of $S^1$-equivariant isotopies, it is equivariant and hence covers an isotopy $\Psi$ of homeomorphisms of $B$. 
Further there exist sequences of Hamiltonian isotopies $\Psi_{G^j}$ and $\Psi_{K^j}$ of $(B,\omega)$ that are generated by smooth Hamiltonian functions $G^j$ and $K^j:B\times [0,1] \to \R$ where $G^j_t \circ \pi = -H^j_t$ and $K^j_t \circ \pi = -F^j_t$ with the following convergence properties.
Because of the split metric $g_M = \pi^*g_B + \alpha \otimes \alpha$ the isotopies $\Psi_{G^j}$ and $\Psi_{K^j}$ converge in $C^0$ to $\Psi$, and because $\| G^j \|_{\rm Hofer} = \int_0^1 \osc(G^j_t)\, dt \leq \|H^j \|$ and $\| K^j \|_{\rm Hofer} = \int_0^1 \osc(K^j_t)\, dt \leq \|F^j \|$ there exist topological Hamiltonians $G$ and $K \in L^1([0,1], C^0(B))$ such that $\|G^j - G \|_{\rm Hofer} \to 0$ and $\|K^j-K\|_{\rm Hofer} \to 0$.

By the results of Viterbo \cite{viterbo:ot06} and Buhovsky-Seyfaddini \cite{buhovsky:ugh10}, the mean-value normalized Hamiltonians $\widetilde G_t = G_t - c(G_t)$ and $\widetilde K_t = K_t - c(K_t)$ are equal for almost every $t \in [0,1]$.
Hence the mean-value normalized contact Hamiltonians $\widetilde H_t = H_t - c(H_t)$ and $\widetilde F_t = F_t - c(F_t)$ satisfy $\widetilde H_t = \widetilde F_t$ for almost every $t \in [0,1]$, since $\widetilde H_t = - \widetilde G_t \circ \pi$, and $\widetilde F_t = - \widetilde K_t \circ \pi$.

Let $\varPhi_R = \{ \varphi^t_R \}$ denote the Reeb flow.
If we denote by $a(t) = - \int_0^t c(H_s)\, ds$ and $b(t) = -\int_0^t c(F_s)\, ds$, then Lemma~\ref{lem:top_reeb} and Theorem~\ref{thm:topgroup} imply that $\{ \phi_t \circ \varphi_R^{a(t)}\}$ and $\{ \phi_t \circ \varphi_R^{b(t)} \}$ are topological strictly contact isotopies associated to $\widetilde H$ and $\widetilde F$, respectively.
Because $\widetilde H = \widetilde F$, Lemma \ref{thm:non-triviality} implies that for all $t \in [0,1]$,
	\[ \phi_t \circ \varphi_R^{a(t)} = \phi_t \circ \varphi_R^{b(t)}. \]
The periodicity of the Reeb flow implies that for all $t \in [0,1]$,
	\[ b(t) - a(t) = \int_0^t c(H_s)\, ds - \int_0^t c(F_s)\, ds \in \Z,\] 
but, because $a^j(t) = -\int_0^t c(H^j_s)\, ds$ and $b^j(t) = -\int_0^t c(F^j_s)\, ds$ converge uniformly to $a(t)$ and $b(t)$, respectively, the function 
	\[t \mapsto \int_0^t c(H_s - F_s)\, ds\] 
is continuous, and thus must be constant.
Therefore $c(H_t) = c(F_t)$ for almost all $t \in [0,1]$, and this implies that $H = F$.
\end{proof}

\bibliographystyle{amsplain}

\def\cprime{$'$}
\providecommand{\bysame}{\leavevmode\hbox to3em{\hrulefill}\thinspace}
\providecommand{\MR}{\relax\ifhmode\unskip\space\fi MR }
\providecommand{\MRhref}[2]{%
  \href{http://www.ams.org/mathscinet-getitem?mr=#1}{#2}
  }
\providecommand{\href}[2]{#2}

\end{document}